
\documentstyle{amsppt}
\magnification1200
\pagewidth{6.5 true in}
\pageheight{9.25 true in}
\NoBlackBoxes

\document

\topmatter
\title
Lower bounds for moments of $L$-functions
\endtitle
\author
Z. Rudnick and K. Soundararajan
\endauthor
\address School of Mathematical Sciences, Tel Aviv University,
Tel Aviv 69978, Israel
\endaddress
\email{rudnick{\@}post.tau.ac.il}
\endemail
\address{Department of Mathematics, University of Michigan, Ann Arbor,
Michigan 48109, USA} \endaddress \email{ksound{\@}umich.edu}
\endemail

\thanks{This research was partially supported by Grant No. 2002088 from the
United States-Israel Binational Science Foundation (BSF),
Jerusalem, Israel. The second author is partially supported by the
National Science Foundation.}
\endthanks

\abstract The moments of central values of families of $L$-functions
have recently attracted much attention and, with the
work of Keating and Snaith, there are now precise
conjectures for their limiting values.  We develop a simple
method to establish lower bounds of the conjectured
order of magnitude for several such families of $L$-functions.
As an example we work out the case of the family of
all Dirichlet $L$-functions to a prime modulus.

{ Resum{\' e}}. Les moments des valeurs centrales des familles
de fonctions $L$, ont suscit{\' e} beaucoup d'int{\^{e}}ret
recemment et, apr{\' e}s le travail de Keating et Snaith,
il y a maintenant des conjectures précises sur leur valeur.  Nous
developpons
une technique simple pour {\' e}tablir des bornes inf{\' e}rieures
sur l'ordre de grandeur
conjectur{\' e} de plusieurs familles de fonctions $L$.
Un exemple est de consid{\' e}rer la famille des fonctions
$L$ de Dirichlet de conducteur
premier.
\endabstract

\endtopmatter
\def\sumstar{\sideset \and^{*} \to \sum}

\noindent A classical question in the theory of the Riemann zeta-function
asks for asymptotics of the moments
$\int_1^T |\zeta(\tfrac 12+it)|^{2k}dt$
where $k$ is a positive integer. A folklore conjecture
states that the $2k$-th moment should
be asymptotic to $C_k T(\log T)^{k^2}$ for a positive constant $C_k$.
Only very recently with the work of Keating and Snaith [4]
modeling the moments of
$\zeta(s)$ by moments of characteristic polynomials of random matrices
has a conjecture emerged for the value of $C_k$.  This conjecture
agrees  with classical results of Hardy and Littlewood, and Ingham (see
[7]) in the cases $k=1$ and $k=2$, but for $k\ge 3$ very little is known.
Ramachandra [6] showed that
$\int_1^T |\zeta(\tfrac 12+it)|^{2k} dt \gg T(\log T)^{k^2}$ for
positive integers $2k$, and Heath-Brown [2] extended this for
any positive rational number $k$.  Titchmarsh (see Theorem 7.19 of [7])
had previously obtained a smooth version of these lower bounds
for positive integers $k$.

Analogously, given a family of $L$-functions an important problem
is to understand the moments of the central values of these
$L$-functions.  Modeling the family of $L$-functions using the
choice of random matrices ensembles  suggested by Katz and Sarnak
[3] in their study of low-lying zeros, Keating and Snaith [5] have
advanced conjectures for such moments. We illustrate these
conjectures by considering three prototypical examples. The family
of Dirichlet $L$-functions $L(s,\chi)$ as $\chi$ varies over
primitive characters $\pmod q$; this is a `unitary' family and it
is conjectured that
$$
\sumstar\Sb \chi \pmod q \endSb |L(\tfrac 12,\chi)|^{2k} \sim
C_{1,k} \ q  (\log q)^{k^2} \tag{1}
$$
where $k \in {\Bbb N}$ and $C_{1,k}$ is a specified positive constant.
The family of quadratic Dirichlet $L$-functions $L(s,\chi_d)$
where $d$ is a fundamental discriminant and $\chi_d$ is the
associated quadratic character; this is a `symplectic'
family and it is conjectured that
$$
\sum_{|d|\le X} L(\tfrac 12, \chi_d)^{k} \sim C_{2,k} X (\log X)^{k(k+1)/2}
\tag{2}
$$
where $k \in {\Bbb N}$ and $C_{2,k}$ is a
specified positive constant.  The family of quadratic
twists of a given newform $f$, $L(s,f \otimes \chi_d)$ (the
$L$-function is normalized so that the central point is $\tfrac 12$);
this is an `orthogonal' family and it is conjectured that
$$
\sum_{|d|\le X} L(\tfrac 12, f\otimes \chi_d)^{k}
\sim C_{3,k} X (\log X)^{k(k-1)/2} \tag{3}
$$
where $k \in {\Bbb N}$, $C_{3,k}$ is a specified constant which depends
on the form $f$.

While asymptotics in (1)-(3) are known for small values of $k$,
for large $k$ these conjectures appear formidable.  Further, the
methods used to obtain lower bounds for moments of $\zeta(s)$
do not appear to generalize to this situation.  In this note
we describe a new and simple method which furnishes lower bounds
of the conjectured order of magnitude for many families of
$L$-functions, including the three prototypical examples given above.
As a rough principle, it seems that whenever one can evaluate the
first moment of a family of $L$-functions (with a little bit to
spare) then one can obtain good lower bounds for all moments.
We illustrate our method by giving a lower
bound (of the conjectured order of magnitude) for
the family of Dirichlet characters modulo a prime.

\proclaim{Theorem}  Let $k$ be a fixed natural number.  Then for
all large primes $q$
$$
\sum\Sb \chi \pmod q \\ \chi \neq \chi_0 \endSb |L(\tfrac 12,\chi)|^{2k}
\gg_k q(\log q)^{k^2}.
$$
\endproclaim



Proofs of the corresponding lower bounds for several other
families, together with other applications of this method, for
instance to fluctuations of matrix elements of Maass wave forms in
the modular domain, will appear elsewhere.  We remark also that we
may take $k$ to be any positive rational number $\ge 1$ in the Theorem.
If $k=r/s (\ge 1)$ is rational, then we achieve this
by taking $A(\chi) = (\sum_{n\le q^{1/(2r)}} d_{1/s}(n)\chi(n)/\sqrt{n})^s$
in the argument below.

\demo{Proof of the Theorem} Let $x:=q^{1/(2k)}$ be a small power of $q$,
and set $A(\chi)=\sum_{n\le x} \chi(n)/\sqrt{n}$.
We will evaluate
$$
S_1 := \sum\Sb \chi \pmod q \\ \chi \neq \chi_0 \endSb
L(\tfrac 12, \chi) A(\chi)^{k-1} \overline{A(\chi)}^k, \qquad
\text{and} \qquad
S_2:= \sum\Sb \chi \pmod q \\ \chi \neq \chi_0 \endSb
|A(\chi)|^{2k},
$$
and show that $S_2\ll q(\log q)^{k^2} \ll S_1$.
The Theorem then follows from H{\"o}lder's inequality:
$$
\sum\Sb \chi \pmod q\\ \chi \neq \chi_0 \endSb  |L(\tfrac 12,\chi)|^{2k}
\ge \frac{|S_1|^{2k}}{S_2^{2k-1}} \gg q(\log q)^{k^2}.
$$

If $\ell \in {\Bbb N}$ then we may write $A(\chi)^{\ell}
= \sum_{n\le x^\ell} \chi(n) d_{\ell}(n;x)/\sqrt{n}$
where $d_{\ell}(n;x)$ denotes the number of ways of
writing $n$ as $a_1 \cdots a_\ell$ with each $a_j \le x$.
As usual $d_\ell(n)$ will denote the $\ell$-th divisor
function, and note that $d_{\ell}(n;x)\le d_\ell(n)$
with equality holding when $n\le x$.

We start with $S_2$.  Note that $A(\chi_0) \ll \sqrt{x}$
and so
$$
\align
\sum\Sb \chi \pmod q\\ \chi \neq \chi_0 \endSb |A(\chi)|^{2k} &=
\sum_{\chi \pmod q}
|A(\chi)|^{2k} +O(x^k) \\
&=\sum\Sb m,n \le x^{k} \endSb \frac{d_k(m,x)d_k(n,x)}{\sqrt{mn}}
\sum_{\chi \pmod q } \chi(m) \overline{\chi}(n)
+O(x^k). \\
\endalign
$$
Since $x^k =\sqrt{q} <q$ the orthogonality relation for characters
$\pmod q$
gives that only the diagonal terms $m =n$ survive.  Thus
$$
S_2 = \phi(q) \sum_{n\le x^k} \frac{d_k(n,x)^2}{n} +O(\sqrt{q}).
$$
Since $d_k(n,x)\leq d_k(n)$ and
$\sum_{n\le y} d_k(n)^2/n \sim c_k (\log y)^{k^2}$ for a
positive  constant $c_k$, we find that
$S_2 \ll q(\log q)^{k^2}$, as claimed.

We now turn to $S_1$.  If Re$(s)>1$ then integration by
parts gives
$$
L(s,\chi) = \sum_{n\le X} \frac{\chi(n)}{n^s}
+ \int_X^{\infty} \frac{1}{y^s} d \Big(\sum_{X<n\le y} \chi(n)\Big)
=\sum_{n\le X} \frac{\chi(n)}{n^s} + s\int_X^{\infty}
\frac{\sum_{X <n\le y} \chi(n)}{y^{s+1}} dy.
$$
Since the numerator of the integrand above is $\ll \sqrt{q}\log q$ by
the P{\' o}lya-Vinogradov inequality (see \S23 of [1])
the above expression furnishes an
analytic continuation of $L(s,\chi)$ to Re$(s)>0$.  Moreover we obtain
that
$$
L(\tfrac 12,\chi) =\sum_{n\le X} \frac{\chi(n)}{\sqrt{n}} +
O\Big(\frac{\sqrt{q}\log q}{\sqrt{X}}\Big).
$$
We choose here $X=q \log^4 q$ and obtain that
$$
S_1 =\sum\Sb \chi \pmod q \\ \chi \neq \chi_0 \endSb \sum_{n\le X}
\frac{\chi(n)}{\sqrt{n}}
A(\chi)^{k-1} \overline{A(\chi)}^k +
O\Big( \frac{1}{\log q}\sum\Sb \chi \pmod q
\\ \chi\neq \chi_0 \endSb |A(\chi)|^{2k-1}\Big).
$$
Since $|A(\chi)|^{2k-1} \le 1+ |A(\chi)|^{2k}$ the error term above is
$\ll (q +S_2)/\log q$.  The main term is
$$
\sum\Sb \chi \pmod q \endSb \sum_{n\le X} \frac{\chi(n)}{\sqrt{n}}
A(\chi)^{k-1}
\overline{A(\chi)}^{k}
+O(\sqrt{X} x^{k-\frac 12}).
$$
Recalling that $x=q^{1/(2k)}$ and using the
orthogonality relation for characters we conclude that
$$
S_1 = \phi(q)\sum\Sb a\le x^{k-1} \endSb \sum\Sb b\le x^k \endSb
\sum\Sb n\le X \\ an\equiv b\pmod q\endSb
\frac{d_{k-1}(a;x)d_k(b;x)}{\sqrt{abn}}  +O\Big(\frac{S_2}{\log q}\Big).
$$

The main term above will arise from the diagonal terms
$an=b$.  Let us first estimate the contribution of
the off-diagonal terms.  Here we may write $an=b+q\ell$ where
$1\le \ell \le Xx^{k-1}/q = x^{k-1}(\log q)^4$.   The contribution
of these off-diagonal terms is
$$
\ll q \sum_{b\le x^k} \frac{d_{k}(b;x)}{\sqrt{b}}
\sum_{\ell \le x^{k-1} (\log q)^4}
\frac{1}{\sqrt{q\ell}} \sum_{an=b+q\ell} d_{k-1}(a;x)
\ll q^{\frac 12 +\epsilon} x^{\frac{k-1}{2} +\frac k2} \ll
\frac{q}{\log q}
$$
since $\sum_{an=b+q\ell} d_{k-1}(a;x)\le d_k(b+q\ell) \ll
(q\ell)^{\epsilon}$.
Therefore
$$
S_1= \phi(q) \sum_{b\le x^k} \frac{d_k(b;x)}{b} \sum\Sb a\le x^{k-1} ,
n\le X \\
an=b \endSb
d_{k-1}(a;x) + O\Big(\frac{S_2}{\log q}\Big).
$$
Since
$$
\sum\Sb a\le x^{k-1}, n\le X\\ an=b\endSb d_{k-1}(a;x)
\ge \sum\Sb a\le x^{k-1}, n\le x\\ an=b\endSb d_{k-1}(a;x) =
d_{k}(b;x),
$$
and $d_k(b;x)= d_k(b)$ for $b\le x$, we deduce that
$$
S_1 \ge \phi(q)  \sum_{b\le x^k} \frac{d_k(b;x)^2}{b} +
O\Big(\frac{S_2}{\log q}\Big) \ge
\phi(q)\sum_{b\le x} \frac{d_k(b)^2}{b} + O\Big(\frac{S_2}{\log x}\Big)
\gg q(\log q)^{k^2}.
$$
This proves the Theorem.

\enddemo

\enddocument